# JORDAN s-IDENTITIES IN THREE VARIABLES

## S. R. Sverchkov*

## INTRODUCTION

Let $J[X_n]$, $SJ[X_n]$, and $As[X_n]$ be a free Jordan algebra, a free special Jordan algebra, and a free associative algebra on a set of generators $X_n = \{aq,\ rzq,\ \ldots,\ x_n\}$ and $S_n$ be the kernel of a canonical homomorphism $7\Gamma : J[X_n] \longrightarrow SJ[X_n]$. Nonzero elements in $S_n$ are called *s-identities*. The Shirshov-McDonald theorem (see [1]) states that if $f(x\ 1,^2,^3)\ G\ S_3$ and $d_{X_3}(f)\ ^\wedge\ 1$ then $/ = 0$. In particular, $S_2 = (0)$. In [2], it was proved that $S_3\ \phi\ (0)$. In [3], 5-identities $Gg,\ Gg\ G\ S_3$ of degrees 8 and 9 were constructed based on the proof in [2]. Currently we know of just a few examples of 5-identities (see [3-7]).

In the present paper we prove that all Jordan 5-identities in three variables are consequences of the Glennie 5-identity $Gg$ of degree 8, which solves K. A. Zhevlakov's problem (see [8, Problem 1.40]). All algebras are treated over a field $F$ of characteristic 0. For standard definitions and notation, we ask the reader to consult [1].

## 1. SHIRSHOV-COHN POLYNOMIALS IN ALGEBRA $J[X_3]$

Let $J = J[X_3]$, $SJ = SJ[X_3]$, $As = As[X_3]$, and $H = Я[X_3]$, where $H[X_3] = \text{ff}(As[X_3], *)$ is an algebra of symmetric elements of $As[X_3]$ with respect to a standard involution $*$. Multiplication operations in the algebras $J$ and $SJ$ are denoted by $\bullet$ and $\circ$. To represent nonassociative words, we use left-normed bracketing.

* Supported by FAPESP (grant No. 2010/50171-8/).

Let $u = yf \ldots .yf$ G As, where $yi$ G $X_3$ and $yi \; \varnothing \; y_{2+i}$, $i = 1, \ldots, k — 1$. The number $h\{u) = \kappa$ is called the *height of a polynomial u*. We assume that $h\{u) = 0$ if $u = 0$. The total number of $x\backslash$ occurring in the representation of $u$ is denoted by $h_{Xl}(u) = e$ and is referred to as the *height of u with respect to x\\*. The beginning and the end of a monomial $u$ are denoted by $b(u) = y\backslash$ and $e\{u) = yk$■ To reduce the number of indices and degrees in representations of polynomials, monomials like $xf$, $x^\wedge$, $xf$ are conventionally denoted by $x_s$, $y_s$, $z_s$, where s, , $*_2$, $*_3$ G N. By convention, for instance, $y_3X\backslash X_2$, $y_2\%szi$, $yiX_2Z\%$ will denote some monomial of the form $xfxfxif$. Below it will be clear that such a representation does not lead to contradictions and simplifies the presentation.

We know from [f] that the set $(\frac{3}{4} = yf \ldots .y'f$, $yi$ G $A_3$, $yi \; \varnothing \; yi+i$, i = 1,…,\frac{3}{4} -1, i G N} of all monomials forms a basis for As and the set $\{\{u\} = u_l + u^*$, $i$ G N} forms one for *SJ*. Fix a basis $u_l$, $i$ G N, for the algebra As. Now, using the Shirshov-Cohn algorithm (see [1]), we construct a set $\{(u_l)$ G J, $i$ G N}, where $m_2(u_l) = \{u_l\}$, $i$ G N. Elements $(u_l)$, $i$ G N, are called *Shirshov-Cohn polynomials*.

A construction algorithm for Shirshov-Cohn polynomials includes a set of rules for each type $u_l$, $i$ G N, and is defined by induction on the height $\kappa(u_l)$, $i$ G N. Moreover, monomials $x_s$, $y_s$, and $z_s$ will have equivalent occurrences in a definition of the algorithm. Therefore, each rule of the algorithm defining some polynomial $\{u_l\}$ extends to all monomials $Uj$ obtained from $u_l$ by permuting $x_s$, $y_s$, and $z_s$.

Below $R_a$ and $U_{a>}b$ are standard operators for Jordan multiplication, i.e.,

$$dRb — U • 6, \; U_a{}^\wedge b — RaRb \; T \; RbRa \; Ra\text{-}bi$$

and $D_ab$ is a Jordan derivation, i.e., $cD_ab = c(R_aRb — RbRa)\text{-}$

### Construction Algorithm for Shirshov-Cohn Polynomials

(1) Let $\kappa(u_l) = 1$; then $(u_l) = 2u_l$.

(2) Let $\kappa(u_l) = 2$ and $u_l = x\backslash y\backslash /$ then $\{x\backslash y\backslash) = 2x\backslash$ ■ $y\backslash$.

(3) Let $\kappa(u_l) = 3$; if $u_l = x\backslash y\backslash Z\backslash$ then $(u_l) = 2yiU_{XltZl'}$ and if $u_l = x\backslash y\backslash X_2$ then $(u_l) = (xm)$ ■ $x_2 + (yix_2)$ ■ $xi - \backslash\{y\backslash x_2xi) - \backslash\{x_2xiyi)$.

Suppose that all Shirshov-Cohn polynomials of height smaller than $n$ are already defined and $\kappa(u_l) = n$.

(4) Let $u_l = X\backslash UjX_2$, where $h(uj) = n — 2$; then

$$<*> = <*_{\langle i \rangle} \bullet^{x_2} + <W> ■ \qquad\qquad - !<w.>.$$

Note that $h(x\backslash X_2 Uj), h(ujX_2Xi) < n$.

(5) Let $u_l = x\backslash Ujy\backslash$ and $h(uj) = n — 2$; then we put $(u_l) = (u^*)$.

(6) Let $u_l = x\backslash Ujy\backslash$ and $h(uj) = n — 2$, with $b\{uf) = y_2$ or $e(uj) = X_2$- Then

$$(v,i) = 2\{uj\}U_{Xl>yi} —$$

$$(yi UjXi).$$

Note that $h(y \setminus UjX \setminus) < n.$

(7) Let $u_{\text{ч}} = xiZiX_2UjZ_2yi$, $h(uj) = n - 5$, and $(uj = 0$ or $e(uj) = x_3)$. Then

$$('Ш) = 2(ziX_2UjZ2yi) \text{ -}Xi \text{ - } (ZiX_2UjZ_2yiXi).$$

The polynomial $(z \setminus X_2 UjZ_2 y \setminus X \setminus)$ was defined in (6).

(8) Let $u_{\text{ч}} = xiZ \setminus X_2 Ujy_2 Z_2 yi$, where $h(uj) = n - 6$. Then

$$('Ш) = \{ZlX_2Ujy_2Z_2yi) \text{ - } Xi + \{XiZiX_2Ujy_jZ_2) \ '2/1 \text{ - } {}^\wedge \qquad ) \text{ - } {}^\wedge(2/1{}^\wedge1{}^\wedge2{}^\wedge2\frac{3}{4}).$$

The polynomials $(z \setminus X_2 Ujy_2 Z_2 y \setminus X \setminus)$ ᵃⁿd $(2/1{}^\wedge1{}^\wedge1{}^\wedge2{}^\wedge2/2\frac{3}{4})$ were defined in (6).

(9) Let $u_{\text{ч}} = x \setminus ZiyiUjX_2Z_2y_2$, where $h(uj) = n - 6$. Then

$$(\blacksquare Ui) = 2(z_1y_1U_jX2Z_2)U_{xum} \text{ - } (y_2ZiyiU_jX_2Z_2X_1).$$

The polynomial $(2/2{}^\wedge1//1{}^{\wedge\wedge}2{}^\wedge2{}^\wedge1)$ was defined in (8).

(10) In defining $(u_{\text{ч}})$, the rule $\kappa$, $\kappa {}^\wedge 9$, will be applied unless $(u_{\text{ч}})$ accepts rules $1, 2, \ldots, k - 1$.

**Examples.** Consider

$$(x_iy_iX_2U_iz_ix_sz_2) = 2(y_ix_2u_iz_ix_3)U_{xuZ2} \text{ - } (z_2yiX_2U_iz_i(x_3x_1)),$$
$$(x_iy_iX_2U_iZ_iy_2Z2) = (yiX_2U_iZ_iy_2Z2) \ \blacksquare \ X_3 + ({}^\wedge1//1{}^\wedge2{}^{\wedge\wedge}1//2) \bullet {}^\wedge2$$
$$\text{ - } {}^\wedge(yix_2u_iz_iy_2z_2x_1) \text{ - } {}^\wedge(z_2x_iy_ix_2u_iz_iy_2).$$

**PROPOSITION 1.** (1) Rules (1)-(10) uniquely define all polynomials $(u_{\text{ч}})$, $i$ G N.

(2) It is true that $\{u_{\text{ч}}) = \{u^*)$ for all $i$ e N.

(3) It is true that $7m((u_{\text{ч}})) = \{u_{\text{ч}}\}$ for all $i$ G N.

(4) The set $(u_{\text{ч}})$, $i$ G N, is linearly independent in J.

**Proof.** Statement (1) follows from the definition of the algorithm; (2) and (3) are easily proved by induction on the height of $u_{\text{ч}}$; (4) is a consequence of (3). The proposition is proved.

In what follows, every relation having occurrences of $(u_{\text{ч}})$, $x_s$, $y_s$, and is a series of relations under all permutations of $x_s$, $y_s$, and $z_s$.

## 2. CONSTRUCTION PRINCIPLE FOR
## THREE-VARIABLE s-IDENTITIES

Let $A$ C J. Denote by $(A)p$ an F-submodule of $J$ generated by $A$, by $(A)j$ an ideal generated by $A$, and by $T(A)$ a T-ideal of $J$ generated by $A$. Put $S = (2(u_{\text{ч}})K_a - (au_{\text{ч}}) - (u_{\text{ч}}a), \ 2(u_{\text{ч}})u_a, \text{ь} - (au_{\text{ч}}b) - (bu_{\text{ч}}a); \ \varnothing$ G N, a, $b$ G $\{x_s, y_s, \ z_s\})F$- Consider / G J. Let $7r(/) = J2^ai\{{}^ui\}$-$_i$

**Definition.** An element $5(f) = / {\text -}^\wedge u_{\text{ч}}(u_{\text{ч}})$ is called an $s\text{-}error$ of /.

Obviously, $7\Gamma(5(f)) = 0$ for every $/$ e J, i.e., $5(f)$ £ $S_3$. Therefore, for any $f$ £ $J$,

$$/ = 5(f) + \sum_i \hat{}ai(m),$$

where $5(f)$ £ $S_3$. The theorem below sheds light on a construction principle for three-variable s-identities.

**THEOREM 1.** We have $+_3 = (S)j$.

**Proof.** Proposition f implies $7T((S)J) = 0$. Therefore, $(S)j$ C $S_3$.

We argue for the inverse inclusion. Let $/$ be a homogeneous element of $S_3$. By induction on the length $d(f)$ of $/$, we state that $/$ £ $(S)j$. The statement is obvious for $d(f) = 2,3$. Let $d(f) > 3$. It is known that the multiplication algebra $R(J)$ is generated by a set $\{R_a, U_{a>b}; a,b$ £ $X_3\}$ of operators (see [I]). Therefore, $/ = \hat{}/Д_{a_i} +$ where $fi$ and $gj$ are homogeneous polynomials of degrees $d(f) — 1$ and $d(f) — 2$, сц, $bj$, $Cj$ £ $X_3$. We have

$$fi = 5(fi) + \sum_\kappa J2a_k(u_k), \quad gj = 5(gf) +$$

where $5(fi), 5(gj)$ £ $S_3$. By the inductive assumption, $5(fi), 5(gj)$ £ $(S)j$. Hence $/ = /0 + Y, a_s(u_s)Ra_s + Y)Pt(ut)U_{bt>Ci}$, where $f_*$ £ $(S)j$, $a_s, b_t, c_t$ £ $X_3$, and $a_s, (3_t$ e $F$. Put

$$/1 = \hat{} \ \hat{}z a_s(z(u_s)R_{as} - (a_s u_s) - (u_s a_s)), s$$

$$/2 = \hat{}_t \ Pt(\hat{}\{ut\}U_{bt}, c_t \sim \{btutct) - (c_iutbt)\} \blacksquare$$

By definition, $/1,/2$ £ $S$. Therefore,

$$/ = /0 + /1 + /2 + \sum_i \hat{}$$

where $(u_l)$ are all different. It follows that $7r(/) = \sum_i J2^a ii^u i\} = 0$ and $\frac{3}{4} = 0$ for all $i$. Thus $f$ £ $(S)j$ and $S_3 = (S)j$. The theorem is proved.

**Remark.** Theorem 1 remains valid for any number of variables. In fact, let $\{ui\}$, $i$ £ **N**, be some basis of $SJ[X_n]$, consisting of homogeneous polynomials, and $(u_l)$, $i$ £ N, be preimages of the basis elements $\{u_l\}$, $i$ £ N, in the algebra $J[X_n\backslash$. The proof of Theorem 1 entails $S_n = (+)jpc_n]$-Unfortunately, a constructive basis is not known even for the algebra $SJ[X_4]$.

## 1. AUXILIARY RESULTS

Denote by $R = ((u_l), i$ e N$)p$ an F-submodule of $J$ generated by all Shirshov-Cohn polynomials. We will write $/ = g$ if $/ — g$ £ $R$. Let $L$ C $(S)j$. By writing $/ = g(L)$ we mean that $f - g$ £ $R + T(L)$. Put

$$(\{+\}+?) — \qquad ) + (u_iu_j), (//,; I_{il\ j}\ j) — (ll, ll\ .j) + (ll, il\ j).$$

**LEMMA 1.** An algebra $J$ satisfies the following implications:

(1)  $/ = 0$, $71\text{-}(/) = o \Rightarrow / = 0$;

*(2)  (Ui)R('Uj) = 0 =r- 2(Ui)R('Uj)* $^=$ *{{*$^u$*j}*$^u$*i) + ('Ujj'Uj});*

(2)  $\{u_\ell\}u_{(U_j)}, (u_k) = 0 \Rightarrow$ ■ $2\{Ui\}U_{(U_j}{}^{\wedge}{}_{U_k)} = \{\{Uj\}Ui\{Uk\}\} + \{\{$"life $\}$        $\{lij\}$);

(3)  $(m)R_{(U_{ij})} = 0(L) \Rightarrow 2(m)R_{\{u_{jj}\}} - (\{Uj\}ui) - (Ui\{uj\})$ e $T(L)$;

(4)  $\{u_\ell\}u(U_j), (u_k) = 0(L) \Rightarrow$ ■ $2\{ui\}U_{(U_j}{}^{\wedge}{}_{U_k)} - \{\{u_\ell\}ui\{u_k\}\} - \{\{uk\}ui\{uj\}\}$ £ $T(L)$.

Proof. (1) Let $/ = J2^a i(^u i) >$ where $(u_\ell)$ are all different. Then $TT(/) = J2^a ii^u i\} = 0\text{-}$

Therefore, $cq = 0$ for all $i$.

(2)  Put $/ = 2\text{-} (\{u_\ell\}u_\ell) - (\text{-}u_\ell\{u_\ell\})$. Then $/ = 0$, $ir(/) = 0$. Therefore, $/ = 0$.

(3)  Is proved similarly to (2).

(4)  Let $2\{ui\}R^{\wedge}{}_u. ^{\wedge} = l + r$, where $l$ £ $T(L)$ and $r$ £ $R$. Then $/ = 0$ and $ir(/) = 0$ for $/ = 2(ui)R_{(U_j)}$
$- (\{u_\ell\}u_\ell) - (ui\{uj\}) - l$. Therefore, $/ = 0$.

(5)  Is proved similarly to (4). The lemma is completed.

Many elements of $S$ are equal to zero in general. Some zero elements of $S$ will be specified in the next lemma. By writing $x\backslash u_\ell$ and $ujx\backslash$ we mean that $b(u_\ell)$ $\Phi$ $x_s$ and $e(uj)$ $\Phi$ $x_s$.

**LEMMA 2.** An algebra $J$ satisfies the following relations:

$$(xiUi)R_{X2} = \sim(u_ix_2)R_{Xl}, \ (zix_2UiX_3z_2yi)R_{Xl} = 0,$$
$$(zix_2u_iy_2z_2yi)R_{xl} = \sim\{xiZix_2Uiy_2z_2)R_{yi}, \ (u\phi u_{x_uy_1} = 0, \tag{1}$$

where $Uj$ $\phi$ $x_2u_\kappa x_3, y_2u_\kappa y_3, zix_2u_kx_3z_2, ziy_2u_ky_3z_2.$

**Proof.** The required result follows immediately from the definition of Shirshov-Cohn polynomials. The lemma is completed.

We will need the following known identities (see [1, 7]):

(2)

(3)

(4)

(5)

(6)

(7)

where $a$, $b$, and c belong to a strongly associative subalgebra of the whole algebra J. Now we introduce the following notation:

$$R(n) = T(2(ui)R_a - (au_\ell) - (u_\ell a), \ \kappa(u_\ell) ^{\wedge} n),$$

$$lJ(n) = T(2(ui)U_{a,b} - (au_\ell b) - (bu_\ell a), \ \kappa(u_\ell) ^{\wedge}$$

where $i$ G N, $a$, $b$ G $\{x_s, y_s, z_s\}$, and n ^ 1. It is easy to see that

$$(ui)R_a = 0 \ (R(n)), \ (Ui)U_{a>b} = 0 \ (U$$

$(n))$, with $h(ui)$ ^ $n$ and a, $b$ G $y_s$, $z_s\}$.

**Definition.** Jordan polynomials of the form $f_n = 2\{zx\backslash y\backslash \ldots x_n y_n)$ ■ $z - \{z^2 X\backslash y\backslash \ldots x_n y_n) -$ $(zxiyi\ldots x_n y_n z)$, $n$ ^ f, are called $k$-*identities*. A T-ideal generated by fc-identities is denoted by $K$, i.e., $K = T(f_n, \ n$ ^ 1).

The definition above implies that $K$ C A₃. Note that fc-identities are identities in all special Jordan algebras.

In what follows, a reference to the Shirshov-McDonald theorem (see [I]) is marked with $(SM)$. Notice that

$$/1 = 2(zx_1 y_1) ■ z - (z^2 x_1 y_1) - \{zxiyiz\} \ 0.$$

Below we show that

$$/2 = 2\{zx_2 y_2 xy\} ■ z - \{z^2 x_2 y_2 xy) - \{zx_2 y_2 xyz\} \ \varPhi \ 0,$$

i.e., is an s-identity; moreover, all s-identities are consequences of $/2$. Define $W(n) = T(R(n), \ U(n - 1), \ K)$, where $n$ ^ 2.

Obviously, $R(4) = 0$ and $[/(3) \ \ = 0$; the s-identity $/2$ shows that $R(5) \ / \ 0$. Theorem I

implies that $S_3 = \overset{\infty}{\underset{n=2}{W(ri)}}$ and $m_2(W(n)) = 0$ for any $n$ e N.

The main property of a T-ideal $W(n)$ is given in

**LEMMA 3.** For any $(u_i)$, $(Uj)$ G $R$ with $\kappa(u_i) = n$ and $h(uj) = m$, we have

$$2(u_i) ■ (uj) = (\{uj\}ui) + (u_i\{u_j\})(\backslash \yen(n + m - 1)). \tag{8}$$

**Proof.** In view of Lemma 1, it suffices to show that $(u_i) ■ (uj) = 0(\text{IT}(n + m - 1))$. We use induction on $m$. If $m = 1,2$ then

$$(u_i)R_a = 0(\text{IT(n)}),$$
$$(ui) ■ (ab) = 2(m)R_{a \cdot b} = 2(11)(¾¾ + R_b R_a - U_{a>b}) = 0(W(n + 1)),$$

where $a, \ b$ e $\{x_s, y_s, z_s\}$, by definition. Suppose that our statement is true for $m - 1$. According to [1], the algebra $R(J)$ is generated by a set $\{R_a, U_{a \cdot b}; \ a, b$ e $\{x_s, y_s, \ z_s\}\}$ of operators, and by the inductive assumption,

$$(■Ui) ■ (Uj) = \{Ui\}R_{\{u_{j\}}} = \text{^} \ a_s(u_s)R_{as} \ Pt\{u_t\}U_{buCt}(W(n \ + \ m - 2)),$$

where $a_s$, $bt$, $Ct$ G $\{x_s, y_s, z_s\}$, $h(u_s) = n + m - 1$, and $\kappa(u_i) = n + rri - 2$. Consequently, $(u_i)$ ■ $(uj) = 0(\text{IT}(n + m - 1))$ by definition. The lemma is completed.

**Definition.** Suppose an s-identity $l$ belongs to $T(R(ri), U(n-1))$, $n \wedge 3$. A polynomial $l$ G $T(R(n-1), U(n-2))$ is called a *reduction of an s-identity* $f$ if $l = l + k$, where $\kappa$ G $K$.

In the next lemma, we construct first instances of reduction. For brevity, Lemma 1 and relation (8) will be used below without further comment.

LEMMA 4. Let $\kappa(u_l) = n$, $3$. Then

$$(m)R_{Xl} = 0(W(n-l)) \tag{9}$$

if (1) $m = x_2 U j X_3$ or (2) $m = y_i U j y_2$;

$$\{u_i)U_{Xl>y_i} = 0(W(n)) \tag{10}$$

if (3) $u_l = x_2 U j X_3$;

$$\{u_i)U_{Xl>X2} = 0(W(n)) \tag{11}$$

if (4) $u_l = x_3 U j X_4$ or (5) $u_l = x_3 U j y \backslash$.

**Proof.** (1) We have

$$(x_2 U j X_3)R_{Xl} = 2(u j X_3)R_{X2}R_{Xl} = (x_i U j X_3)R_{X2}(R(n-1)).$$

Consequently, the polynomial $(x_2 U j X_3)R_{xl}$ is symmetric in $x\backslash, x_2, x_3$ modulo $R + R(n-1)$. Now

$$\begin{aligned}
\{u_j x_3)R_{Xl}R_{X2} &= 2 (u_j) \ R_{X3}R_{Xl}R_{X2} \ \{x_3 U j')R_{Xl} \ R_{X2} \\
&= \{x_3 U j')R_{Xl} \ R_{X2} = \{x_3 U j X\backslash)R_{X2} \\
&= ---(x_2 U j X_3)R_{xl}(R(n-1))
\end{aligned}$$

and $(x_2 U j X_3)R_{Xl} = 0(R(n-1))$.

(4), (5) The relation proved above readily implies that

$$(x_i U j X_2)U_{X3>X4} = 0(R(n-1)), \tag{12}$$

where $h\{x\backslash U j X_2) = n$, and

$$\{x_3 U j y_i)U_{X2tXs} = 0(R(n-1)), \tag{13}$$

where $h\{x_3 u j y_i) = n$.

(2) Let $x\backslash = x^e$, $x_2 = x^k$, and $x_3 = x^m$, with $e, k, m$ G N. Assume m $\wedge$ k. If $\kappa = m$, then

$$(x^k U j X^k)U_x e_{tyi} = \{u_j) U_x k \ _x k U_x e_{>yi} = 2(u_j)R_x k U_x k + e_{yi} = 0(R(n-2), U(n-2)).$$

If $m > k$, $m = m_i + k$, then, similarly,

$$\{x^k U j X^{mi+k})l J_x z_{tyi} = \{u_j X^{mi})U_x k \ _x k U_x e_{>yi} = 2(u_j X^{mi})R_x k U_x k + e_{yi}$$

$$= (x^k U j X^{mi}) U_x k +_{\pounds t y i} (R(n-1), U(n-2)).$$

If $m$ is divisible by $\kappa$, then $(x^k U j X^m) U_x <-\hat{\ }_{yi} = 0(R(n-1), lJ(n-2))$. If $m$ is not divisible by $k$, then (using the Euclidean algorithm for $m$ and $\kappa$) we define $m = \kappa \blacksquare qo + r\backslash, \kappa = r\backslash \blacksquare q\backslash + \Gamma_2, \dots,$ $r_s\text{-}1 = r_s \blacksquare q_s$, where $r_s = (\mathrm{m}, \kappa)$ is the greatest common divisor of $m$ and $k$. By the above,

$$(,X^k U j X^m) U_x e_m = (x^k U j X^{ri}) U_x e i_{tyi} = \{x^{r_2} U j X^{ri}) U_x {}^{\wedge}2_{tyi}$$

$$= \dots = (x^{rs-l} U j X^{rs}) U_{x^.} e_{a>yl} = 0(R(n-1), U(n-2)).$$

(2) Let $y i = y^k$, $y_2 = y^m$ and $m \wedge k$, $m = \kappa + m\backslash$. If $m\backslash = 0$, then

$$(y^k U j y^k) R_x \text{ i} = \{uj) U y k_t y k R_{xl} \{U\{n\text{-}2)).$$

Therefore,

$$(\text{У Т?У}) R_{X1} {}^= \{Uj) U_y k_{\ y} y k \ R_{\chi} 1 \ (Uj)(2 R_y k U_y k_{>\text{Ж}1} \ U_y 2 k^{\wedge}_{xl})$$

$$- (\text{У }^{\wedge} j) R y k_{,\chi} 1 \ \text{T} (U j y) U_y k^{\wedge}_{xl} \ (Uj) U_y 2 k_{>\text{Ж}1}$$

$$= 0(R(n\text{-}2), U(n\text{-}2)).$$

If mi > 0, then

$$(12)$$

and

$$(44^{\wedge}1)\tfrac{3}{4} = (^{\wedge}1)^{\wedge}\tfrac{3}{4} = (^{\wedge})(2\tfrac{3}{4}.^{\wedge} - [/_{y2} y_{\text{Ж}1})$$

$$= 2(^{\wedge})\tfrac{3}{4}.^{\wedge} = \{y^k U j y^{mi}) U_y k_{tXl} + (^{\wedge}44,^{\wedge} = (\text{уЧ-у}^{\text{ТМ!}})!/^{\wedge} ,=$$

$$0(Д(\text{п - 1}), U(n-2)).$$

Thus

$$\{y i U j y_2) R x_l = 0(R(n-1), U(n- \qquad\qquad\textbf{(14)}$$

where $H\{y\backslash u, y_2\}$ 2)), $n$. The lemma is completed.

# 1. /-IDENTITIES

Consider some set of monomials $u\{$ e As, г e /, where $h_z(ui) = 1$. On the Shirshov-McDonald theorem (see [1]), the following representation is well defined:

$$\text{djtlj j} \longrightarrow \text{\underline{'} } OLj$$
$$(Uj). \quad i \qquad i$$

An elegant form to represent the Glennie s-identity of degree 8 was found by I. P. Shestakov (see [7]): namely,

$$Sh = Sh(x, \ y, z) \ = \ ([,r^2, [x, \ y]^3]) - 2([z, \ [x, \ y]^3]) \bullet 2,$$

where $[x, \ y] \ = \ xy — yx$.

Let $u_i$ G As, $i$ G /, be a monomial. Define $[u_i] \ = \ u_i — u^*$. Obviously, $[u_i]$ belongs to $L(\text{As}, *)$, the Lie algebra of skew-symmetric elements of As (see [9]).

**Definition.** Jordan polynomials of the form

$$g_a \ = \ \{[z^2, [\mathbf{a}]]) - 2([z, [\mathbf{a}]]) \bullet z,$$

where $a$ is an arbitrary monomial in $\text{As}[x,y]$, are called *commutator s-identities,* or, for brevity, Sh *-identities*.

In [9], it was proved that all commutator s-identities follow from $Sh(x,y,z)$, and that the s-identities $Gs$ and Sh are equivalent.

**THEOREM 2.** All /-identities follow from $Sh(x,y,z)$.

**Proof.** It suffices to show that /-identities are consequences of Sh-identities. Let $u_i$ be an arbitrary polynomial in $\text{As}[x,y]$. By the definition of Shirshov-Cohn polynomials, we have

$$(zUiZ) \ = \ \{ZUi) \ \blacksquare \ Z + (UiZ) \ \blacksquare \ z - {}^{\wedge}(z^2Ui) - {}^{\wedge}(UiZ^2).$$

Therefore,

$$2(ZUi) \ \blacksquare \ z - (z^2Ui) - (UiZ) \ = \ (ZUi) \ \blacksquare \ Z - (UiZ) \ \blacksquare \ z - {}^{\wedge}(Z^2Ui) + {}^{\wedge}(UiZ^2)$$

$$= \ ((ZUi) - (ZU^*)) \blacksquare z - i((z^2Ui) - (z^2u^*)).$$

On the other hand,

$$2((ZUi) - \{ZU^*)) \ (ZUi) + (v^*z) - (zv^*) - (mz) \ = \ (z[m]) - ([m]z) \ = \ ([z, [m]\}).$$

Consequently,

$$/ = \backslash [[z, \text{M}]) \bullet {}^{\wedge} - {}^{\wedge}([-\text{r}^2, \text{N}]) = \sim \backslash \textit{зиц-}$$

The theorem is proved.

## 5. REDUCTION OF s-IDENTITIES

Here we prove that all s-identities follow from $Sh(a;, y, z)$ and an infinite series of s-identities like

$$gi \ = \ {}_2(X_l \text{У}_l \text{Ш}X_2 Y_2) \ \blacksquare \ Zi - (ZlXiyiUiX_2y_2) \sim$$

$$(xiyiUiX_2 y_2i), \ g_2 = 2(y_1 e_2 u_i x_1 y_2 x_2) \ \blacksquare \ zi - (e_1 y_1 e_2 u_i x_1 y_2 x_2) -$$

$$\partial_3 = {}_2(xiyiu_i{}^,y,X_2)Uy_{3|Z1} - (y_3x \backslash y \backslash u_iy_ix,_2x) - (z \backslash X_1 Y_1 u_iy_ix,y_3) \blacksquare$$

We need the following:

**LEMMA 5.** For any $u_j$, $i$ **G** N, with $1\partial(u_j) = n$, we have

$$(ui)R_{X1} = 0(W(n-1)), \text{ where } u_j = x_2 Ujy \backslash \text{ and } n \wedge 3; \qquad (15)$$

$$(u_j)u_{x1} \cdot y_u y_{x2} = 0(W(n+1)), \text{ where } n \wedge 1; \qquad (16)$$

$$\{ui)U_{x1} \cdot {}_{z1y1X2 \cdot yi} = 0(W\{n+2)), \text{ where } n \wedge 1; \qquad (17)$$

$$(z_1 x_1 z_2 u_i y_1 x_2) \blacksquare y_2 = \sim(zix_1 z_2 u_i y_i x_1) \blacksquare y_2 \{W(n+4)), \qquad (18)$$

where $u_j$ may be missing;

$$(ZiXiZ_2 UiX2y_2) \blacksquare X_3 = (W(n+4)), \qquad (19)$$

where $u_j$ may be missing.

**Proof.** (15) We have

$$(x_2 Ujy \backslash)R_{X1} - (ujy \backslash X \backslash)R_{x_2} - {}_2(ujyi)R_{X1}R_{x_2} \quad (x \backslash Ujy \backslash)R_{X2}.$$

Hence the element $(x_2 Ujy \backslash)R_{X1}$ is symmetric in $x \backslash$, $x_2$ modulo $W(n-1)$, and

$$(x_2 u_j y_1)R_{X1} = -(ujyix_2)R_{X1}\{W\{n-1)).$$

Therefore,

$$\{uiyi)R_{x1}R_{X2} = 0(W(n-1)), (xiUjyi)R_{X2} = (x_2 Ujyi)R_{X2}(W(n-1)). \qquad (20)$$

Thus the polynomial $(x \backslash Ujyi)R_{X2}$ is symmetric in $x \backslash$, $x_2$ modulo $W(n-1) + R$. Below we assume that $x \backslash = x_2$ for brevity. Consider all possible types of $uj$.

(a) Let $Uj = y_2 Uk$ or $Uj = UkX_3$, where $h(uk) = n-3$. Then

$$(x \backslash y_2 Ukyi)R_{X1} = {}^'2(y2Uk)U_{x1t}y_1 R_{x1} \underset{(5)}{=} (y_2 Uk) (U_x{}^\wedge_{yi} + R_{yi}U_{x1jx1}) \underset{(14),(20)}{=} 0\{W(n-1)).$$

Similarly,

$$(x_2 u_k X_3 yi)R_{x1} = 0(W(n-1)).$$

(b) Let $Uj = z \backslash UkZ_2$ or $Uj = z \backslash$. Then

$$(x_1 z_1 u_k z_2 yi)R_{x1} = -(z \backslash UkZ_2 y_2 x \backslash)R_{x1}$$

$$= -2(u_k z_2 y_2)U_{Z1' X1}R_{X1} + (xiu,_k z_2 yiZi)R_{X1} =$$

$$-2\{u_k z_2 yi)(U_x{}_2 + R_{Zi}U_{XuXi}) + (xiu,_k z_2 yiZi)R_{X1}$$
$$\underset{(5)}{}$$

$$\underset{(14),(20)}{=} (xiu_k z_2 yiZi)R_{ui}\{W\{n-1)).$$

If $b(u_k) = xs$ or $l(uk) = zs$, then $(x\backslash Ujyi)$ $R_{X1} = 0(W(n — 1))$ by the argument in (a). By a similar argument applied to $(x\backslash u_k Z2yiZ\backslash)$, either we are led to relation (15) or make the conclusion that the monomial $ujy\backslash$ has the form $z_m y_m \ldots z\backslash y\backslash$. Consequently,

$$2(xiz_m y_m \ldots .ziyi)R_{X1} - \{x\backslash z_m y_m \ldots ziyi) - \{xiz_m y_m \ldots zm) \text{ e } K \text{ C } W\{n — 1).$$

(16)   We have

$$(ui)U_{Xl} \cdot y_{l \ X2} — \{tii)(R_{Xl}U_{yl \ X}2 \text{ T } Ryi U_{Xl \ X2} \ U_{Xi \ yi}R_{X2}) = \qquad 0(W(n \text{ T } 1)).$$
$$\underset{(4)}{} \qquad \underset{(1),(11),(20),(15)}{}$$

(17)   Consider all possible types of $u_l$: namely, (a) $x_3 UjX_4$, (b) $y_2 UjX_3$, (c) yrrc/Уз, and (d) $y_2 UjZ_2$, where $uj$ may be missing.

Type (a). We have

$$\{X_3 UjX_4)U_{x_i} \\ -Zl,X_2-yi \underset{\ulcorner A \urcorner}{—} \{x_j UjX_4)(R_{Xl} \overset{Uzi,X_2 yi}{} \text{ T } \overset{\mathbb{H}}{Z} i^{\wedge}x^i,_{X2}yi^{\cdot} \qquad x_r z_j/-\cdots X_z-yiJ \\ = \ 0(W \quad (n \ + \\ 2)).$$

Type (b). Put $y_2 Uj = u_k$ ■ Then $h(u_k) = (n — 1)$, and

$$4(u/_c X3)Uy_{l \cdot X2_l}z_{l \atop \underset{(4),(8)}{—}} \ \{yiX_z U_k X_z Zi) + (x_z yiu_k x_z zi) + (ZiU_k X3yiX_z)\{W(n + 2)).$$

Therefore,

$$(yix_2 u_k X3Zi)R_{x1} = A(u_k X3)Uy_{l \cdot X2t Z1}Rx_l(^w(n + 2)).$$

Also

$$\overset{\prime\wedge\{U_k X3)Uy_{l \cdot x_2}^\wedge zi Rxi}{\underset{(5)}{\text{------------}}} 4(11^\wedge.^\wedge3)( \ Uy_{r_x 2,XlRzi} + Ryi-X2^z Zl,Xl \ Uyi-_{x_z}Zl-Xl) \\ = \qquad \wedge\{UkX3)Uy_{l \cdot X2l Zl}-x_l(W \ (n+2)).$$

(21)

On the other hand,

$$Ц(u_k x_3)u_{yI} \ ■^{X2,z\backslash -x\backslash} \underset{(5)}{—} \wedge\wedge\backslash U^j kX3)(JJy_{tr}z\backslash -x\backslash Rx_2 \text{ T } U_{x_z}zi-xi \ Ryi \ Rz\backslash-x\backslash Uy_{tr}x_2) \\ \underset{(9),(1),(10),(15)}{} \qquad \{ziXiu_k x_3 yi)R_{X2}(W(n + 2)).$$

Hence

$$\{y\backslash X_2 U_k X_3 Zi)R_{xl} \quad = \qquad (z_i X_i U_k X3yi)Rx_z \qquad = \\ \{y\backslash X_3 U^*_k XiZi)R_{x_z} \qquad = \\ \{zix_z v^*_k xiyi)R_{x_z} = \{yiXiu_k x_z zi)R_{x_3} =$$

If we interchange $2/1$ and $z\backslash$ we have

$$\{y\backslash X_2{}^\wedge kX_3 Z\backslash)\ R_{\chi^i} - \{X_2 UkX^r s)Uy_{l_t Z_l}R_{\chi^2}\qquad i,Z\backslash X_2 UkX_\backslash y\backslash)\ Rxi$$
$$\overline{Ry_s Ux;_{s_t}z}^2 i\ \overset{\{X_2{}'U'_tX_s\}}{\nabla V}(\ Uxi_\bullet z\backslash Ryi\ \ Uy_{l_t}x;iRzi\ +\ Rx_\bullet iUy_{l_t z_t}\ +$$
$$+\ RziUyi_txi)\ \ \{zlX_2 UkX^\wedge y\backslash)R_\chi i$$
$$=-(z_1 x_2 u_k x_3 y_t)(W(n\ +\ 2)).$$

This implies

$$(yix_3 UkX_2 z_1)R_{xl}=-(zix_2 UkX_3 yi)R_{xl}=-(yix_1 u_k x_3 z_1)R_{X2}(W(n\ +\ 2))$$

and

$$(yix_3 u_k x_2 zi)\ R_{Xl}=-(z_1 x_3 u_k x_2 yi)Rx_t(W(n+\ 2)),$$
$$(mx^\wedge+UkX^\wedge z^\wedge Rx^\wedge=(\sim l)^a\ (y_i x_t UkX_2 z_1)R_{X3}(W(n\ +\ 2)) \tag{22}$$

for all $a\ G\ S_3$.

Now if we use relation (22) we obtain

$$a\ =\ \{yix_2 u_k x_3 z_1)Rx_t\underset{(14),(9)}{=}\quad 4((\text{yi o }x_2)u_k(z_t\text{ o }x_3))R_{Xl}$$
$$\underset{(22)}{=}2(((2/1\circ X_2)u_k(z_t\ \circ\ X_3))\ +\ \{\{zi\ \circ\ X_3)u_k\{yi\ \bullet X_2)))R_{Xl}$$
$$=\ 4\{uk)Uyi\ \text{-}X2.ZI\text{-}X3\ R_{Xl}(W(n\ +\ 2)). \tag{23}$$

On the other hand,

$$a\ =\ \{yiX_2 U_k x_3 zi)R_{Xl}\underset{(15)}{=}2\{(yiox_2)u_k x_3 zi)R_{Xl}$$
$$\underset{(15)}{=}\ '^\wedge\{\ 'U^j kX_3)Uy_{t\bullet}X^2{}_,ziRxi\underset{(21)}{=}\ '^\wedge\{/'U'kX_3)Uy_{t\bullet}X^2{}_,zi\text{-}xi\ (\text{W(w + 2)}).$$

Similarly,

$$a--4;\{\$2'U'k')Uy_{t\bullet}xi_tZl'X3$$
$$^\wedge\{\$3^\wedge k)Uv_{t\bullet}.X2_,Zl'Xl(W'(P'\qquad\qquad 2))^*$$

Therefore,

$$2\ a=4(\{\ae_3\ u_k)+(u_k xs))U_{yi}\bullet X_2 yZi\text{-}Xi-8\ (v,k)\ Rx^\wedge Uy\text{-}^\wedge\ \text{-}X2yZi\text{-}Xi$$
$$\underset{(22)}{=}{}_8\{tik)RxiU_{yi\bullet}x^2{}_,zi\text{-}xs(^\wedge{}^\wedge 2).$$

Consequently,

$$2\&\underset{(23)}{=}4(\%)\ [\ U_{yi}\ _{\bullet x2_,zi\text{-}x_,)}\ Rx_{\ 2]}\ -\underset{(3),(22)}{=}\text{-}8\ \{u_k)\ [\ U_{Xl}\ _{,z2\text{-}x_3\ 1}\ Ry\ _{1\bullet\ae_,]}$$
$$\underset{(16)}{=}\frac{\&}{+}\{ll>k)Uxi_{,,}\textbf{21-2:3}Ry\backslash\text{-}X2\ (\text{IC}(n\ +\ 2)),$$

where $[U, R] = UR - RU$. Thus

$$a = L[u_k)u_{x_2} \quad {}_{yZl - X^\wedge Ryi - Xl}$$

(22)

$$= 2(x_2 u_k(z_1 \circ x_3)) + ((z_i \circ x_3)u_k x_2)(U_{yuXl} \cdot R_{Xl}R_{yi} - R_{yi}R_{Xl})$$

(1) (9) TO (15) $\sim \backslash^\wedge yiX2UkX3Zl^\wedge + (^Zl^3{}^ukX2m))Rx_1$

$$\underset{(22)}{=} -a(W(n + 2)).$$

Hence

$$a = (yix_2 u_k x_3 zi)R_{xl} = 0 \ (W \ (n + 2)). \tag{24}$$

Now (21) entails $\{u_k x_3)U_{yi \cdot X2 \cdot Zl \cdot Xl} = 0(W(n + 2))$ and (24) gives rise to (19).

Type (c). We have

$$(y_2 Uj'y_3)U_{Xl \cdot Zlt X3} \cdot y_1 = \{y_2 Ujy_3\}\{U_{Xl \cdot Zl > X_2}Ry_1 \ - \backslash -U_{Xl \cdot Zl > yi}R_{x2} - R_{Xl \cdot Zl}U_{X2 > yi})$$

(5)

$$\underset{(i)}{=} 0(W \ (n + 2)).$$

Type (d). We have

$$'^\wedge\{y2'U'jZ_2'\}U_{Xl \cdot Zljx2yi} \cdot y_1 \ -'^\wedge\{y2'U'jZ_2')(U_{Xl \cdot ZljX2}Ry_1, U_{Xl} \cdot _{Zl} R_{X2} \quad Rxi-z-$$

$$jU_{X2l \cdot y_l}) \ (5)$$

(9) $(¾ (1)$ $+ (^Zl^Xiy2UjZ_2X_2))Ry_1$

$$+ \{ziXiy_2UjZ_2yi)R_{X2} - (ziXiy_2UjZ_2)U_{x2lyi}(W(n + 2)).$$

However,

$$)z\backslash X\backslash y_2 Uj \ z_2')U_{X2l}y_1 - 2\{x\backslash y_2 Uj \ z_2) \ R_{Zl} \ U_{X2l}y_l$$

$$'2'\{X\backslash y_2 UjZ_2)(^\wedge \ Rx_2^\wedge zi, yi \ RyiUx_2, zi + Ux_2, yiRzi$$

(3)

$$+ \ U_{zljy_1}R_{X2} + U_{X2lzi}R_{yi})$$

$$= \{yixziy_2 UjZ_2x_3)R_{zl} + (ziXiy_2 UjZ_2yi)R_{X2}$$

$$+ \ (2^\wedge 1¾¾¾¾)¾^\wedge + 2)).$$

On the other hand,

$$(Ш\text{Ж}2 - 1 \ y2UjZ_2X_2)R_{zl} = \sim\{Xiy_2UjZ_2X_2Zi)R_{yi}.$$

Therefore,

$$4y2UjZ_2)U_{Xl \cdot Zux2 \cdot yi} = ((x_3y_2UjZ_2zi) + (\varepsilon1^\wedge2^\wedge-^\wedge2^\wedge1))^\wedge(^\wedge + 2)). \tag{25}$$

Now we multiply the operator $U_{xi \cdot zi > x2 \cdot yi}$ out in the opposite direction via (4). Similar computations show that

$$4y2UjZ_2)U_{xl} \quad _{-Z_r x_2 \cdot yi} - 4(2/2^{\wedge} - 22)(\frac{3}{4} U_{xi} - zi, X2 N Rx2^\wedge Xl - Zl, yi U_{x 2, yiRxi'Zl})$$

$$= \quad \sim ((X2y2UjZ_2X_1Z_1) + (x_1y2U_jZ2X_2Z_1))Ry_1(W(n+2)).$$
$$\scriptstyle (1),(3),(5),(9),(15)$$

In view of (25), we have $\{y_2 UjZ_2)U_{Xl} \cdot z_1, {}_{X2} \cdot y_1 = O(W\{n+2))$.

(16)  Follows from (25) and (17). The lemma is completed.

Theorem 2 and Lemmas 2, 4, and 5 can be combined to yield

**THEOREM 3.** All s-identities in three variables follow from Sh and an infinite series of s-identities like $gi$, $g_2$, and #3.

Now we prove that all three-variable s-identities of height at most 6 are consequences of Sh. By virtue of the fact that $R(4) = [/(3) {}_{(SM)} = 0 {}_{(SM)}$, there are no s-identities of height 5. It suffices to show that

$$(xiyizix_2y2) \blacksquare z_2 = 0(\text{Sh}), \tag{26}$$

$$(yiz_1x_1y_2x_2) \blacksquare z_2 = 0(\text{Sh}). \tag{27}$$

We have

$$(x_1y_1z_1x_2y2) \blacksquare z_2 {}_{(14)} = 4(zix_2)U_{Xl \cdot yity2}R_{Z2} - (y_2ziX_2yiXi)R_{Z2}(\text{Sh}).$$

On the other hand,

$$4:\{zix_2)U_{xl} \quad _{'yl > y2 R22 \ldots -7. - \ 4(2^{\wedge}2)(- \ U_{x_l'y_{i_1}Z2Ry2} + {}^{\wedge}1-91^{\wedge}242)}$$
$$\scriptstyle (5),(17)$$
$$= -\{xiyiZix_1z_2)R_{y_2} + (xiyiZix_2)U_{y_2 z2}$$
$$\scriptstyle (!),(15)$$
$$= \{y_2XiyiZix_2)Rz_2 \ (\text{Sh}).$$

Therefore,

$$\{xiyiZiX_2y_2) \blacksquare Z_2 = \{y_2X\backslash yiZiX_2) \bullet \underline{2}_2 - \{y_1ZlX_2yiXi) \bullet \underline{2}_2(\text{Sh}).$$

Consequently, (26) follows from (27). Below we need the following identities:

$$aUb, bRt = 2 \ blJ_a, tRb — tUb, bRa, \tag{28}$$

$$aDt, bU_{C)C} = (iD_{tUcct}b + 2aD_{cUtbta}. \tag{29}$$

In view of the Shirshov-McDonald theorem (see [1]), (28) need only be verified in a free associative algebra. We have

$$(babt + tbab) = 2(abt + tba) \ ob — (abtb + btba).$$

Identity (29) follows immediately from a known identity in [1]. We write

$$dDb.c^2 — 2 \ (iDb-c, c-$$

Indeed

$$2\&Dt,(b-c)-c\text{й-}Di\ b-c^2 — 2o{>}Dt\text{-}_{,\text{-}\varphi}b\ \top\ 2Q,Df._\text{c}.i)\ _c\ 2aDt\text{-}(b-c)\ ,_c$$

$$®\text{-}^t\text{-}c^2,6\ 2\ ciDt\text{-}b\text{-}_{\varphi'_c} =$$

$$aD_{tUccb} — 2$$

$$aD_oj j_{bt>c}\text{-}$$

Therefore, Theorem 2 gives rise to

$$(\text{з·}^\bullet\text{Ж}^\bullet\text{·}\text{Y}^\bullet\text{·}\text{Ж}^\bullet\text{·}\text{Y}_2) \cdot зз = 4yiγ \backslash JU \qquad R \qquad (X = y, y, y)$$

Furthermore,

$$— 4\{^\wedge\text{i}\ *\ _{oci)}U^\wedge y_{1X_2}{}^\wedge{}_{Z_2}Ry_2\ ^\wedge y2U^\wedge y_{1X_2}{}^\wedge{}_{Z_2}Rz_1{}'Xi\ (28)$$

$$^2Uz\backslash\text{-}x\backslash^\wedge y2R\{y_jx2)\ 2(^12/1^2)4(1{:}1,2/2\text{-}^2{:}2$$

$$= (\{//\text{l}^\wedge\text{H}^\wedge\text{l}\}^\wedge)\bullet y2 + 2((2/1^2)2/2^2)(^1^3{:}1 + R_{X_i}R_{Z_1}\ (28)$$

$$¢4^1,3{:}1)\ 2((^1^1)^2//2)(\text{-}^\wedge^1\text{-}^3{:}2\ 4\cdot Rx2Ryi\ ^\wedge V_1,X_2{\sim})\ 2XlU(_{Z_1}y_{1x_2})_{iz}2\ Ry_2$$

$$2//2^\wedge4(_ziyiX2){>}_z2^\wedge{}^{1z}i\ 4\cdot\ 2^\wedge2^3{:}1,2/2\text{-}^\wedge(2{:}12/13{:}2)$$

$$= (^22/1^1^1^2)\bullet 2/2 + (2/1^22/2^1)\bullet ^1$$

$$(1),(9),(15)$$

$$- (^1^1^2//2^2)\bullet 2/1 - (^22/1^1^1^2)\bullet 2/2\ 4\cdot 2^\wedge2^3{:}1,2/2\text{-}^\wedge(212/13{:}2)\ (^{\wedge\wedge})\ *$$

Therefore,

$$2Z_2U_{Xl},y_2\ ^\wedge(ziyiX_2)\qquad 2\text{yi}\ \text{C}/^\wedge1\ ,\text{Ж}2\ \text{-}R(XlZ_2X_2)$$

$$- \text{'}2ziU(_{x1Z2}y_2{}^\wedge{}_{yi}R_{X_2}\ \text{-}j\text{-}2x_2U(_{x1Z2}y_2)_{tyi}R_{Z_1}\ (_2H)$$

$$2(xiz_2y_2)U_{zl\}X2}\ R_{yi}$$

$$= (xiz_2y_2x_2yi)\ \blacksquare\ Zi\ \text{-}\ (ziXiz_2y_2x_2)\ \blacksquare\ 2/\text{i(Sh)}.$$

$$(9),(15)$$

$$(,ZiXiyiX_2y_2)\ \blacksquare\ Z_2 = 2((xiz_2y_2x_2yi)\ \blacksquare\ Zi\ \text{-}\ (ziXiz_2y_2x_2)\ \blacksquare\ y\ 1 =$$

$$4(:TCiZ_2y_2x_2yi)\ \blacksquare\ Zi(Sh).$$

On the other hand, $(z\backslash Xiyix_2y_2)\ \blacksquare\ z_2 = 0(Sh)$. Therefore, $(x\backslash z_2y_2x_2yi)\text{-}z\backslash = 0(Sh)$. Consequently,

$$(15)$$

relation (27) and hence (26) will be valid.

## 6. A-REPRESENTATION FOR s-IDENTITIES

Let $/ = f(x,y,z,t)$ G $J[x,y,z\setminus$ and $d_z(f) = 2$. Linearizing $/$ with respect to $z$ yields $I(f) = f(x,y,z,t) + f(x,y,t,z)$. Assume $I(f) = tA(x,y,z)$, where $A(x,y,z)$ G $R(J[x,y, z])$; here $R(J)$ is an algebra of right multiplications for $J$. An expression like $I(f) = tA(x,y,z)$ is called an *R-representation* of the polynomial $/$.

By definition, if $I(f) = tA(x,y,z)$ is an R-representation of $/$, then $g(x,y, z)A(x,y, z)$ G T$(/)$ for any $g(x, y, z)$ G $J[x, y, z\setminus$. If we write (28) and (29) in the A-form and express an element $aЩ∂$ via $bD_a j$ and $tR_a Rb$ we obtain

$$aU_{b>}bRt — 2bU_{a},{}_tR_b \; tUb,bRa,y$$
$$aDt,bU_c,c — ttDfjj_{c\,c}, \; b \; A \tag{30}$$
$$oXJb,t — tR_bR_a \; A \; bD_a \; j.$$

By using identities (30), we can readily construct A-representations for any Jordan polynomials.

Elements of A(5) of height 6 are written thus:

$$g1 = {}_2(zx_iy_ix_2y_2) \; ■ \; z - (z^2x_iy_ix_2y_2) - (zxiyix_2y_2z),$$
$$92 = {}_2(y_2x_iy_izx_2) \; ■ \; z - (zy_2xiyizx_2) - (y_2xiyizx_2z),$$
$$g_3 = {}_2(x_iy_2zy_ix_2) \; ■ \; z - (zxiyizy_ix_2) - (xiyizy_ix_2z),$$
$$g_4 = {}_2(x_2zx_2y_iz) \; ■ \; y_2 - (y_2x\setminus zx_2y\setminus z) -$$

We construct A-representations for the polynomials above.

**PROPOSITION 2.** For any $t = \{u_i\}$ and any $I\varepsilon(u_i) = n) 1$, the following relations hold: $(zty_2x_2yi)$

$$■ \; xi \; A \; (zy_2txiv_i) \; ■ \; x_2 - (zxitv_2x_2) \; ■ \; vi$$
$$z \tag{31}$$
$$A \; (txiyix_2)Uy_2,z - (ty_2x_2yi)U_{XuZ} = 0(W(n \; A \; 3))$$

$$((ztyiXiy_2) - (zty_2xiyi)) \; ■ \; x_2 \; A \; (zx_2ty_2xi) \; ■ \; yi$$
$$- (zx_2tyixi) \; ■ \; y_2 \; A \; ((x_2tyiXiy_2) - (x_2ty_2xiyi)) \; ■ \; z \tag{32}$$
$$A \; ((ty_2xiyi) - (tyiXiy_2))U_{X2>z} = 0(W(n \; A \; 3)),$$

$$((ztyix_2y_2) \; A \; (zy_2tx_2yi)) \; ■ \; Xi - ((zty_2xiyi) \; A \; (zyitxiy_2)) \; ■ \; x_2$$
$$A \; ((y_2x_2tyixi) \; A \; (xityix_2y_2)$$
$$- (y\setminus X\setminus ty_2x_2) - (x_2tyiXiy_2)) \; ■ \; z \tag{33}$$
$$A \; (ty_2xiyi)Ux_2,z - (tyix_2y_2)U_{XuZ} = 0(W(n \; A \; 3))$$

$$\left(2\langle y_2ty_1x_2z\rangle + \frac{1}{2}\langle y_2tx_2y_1z\rangle\right) \cdot x_1$$
$$+ \left(2\langle x_1tx_2y_1z\rangle + \frac{1}{2}\langle x_1ty_1x_2z\rangle\right) \cdot y_2 \equiv 0(W(n+3)) \tag{34}$$

$$((zty,xiyi) + (zy,txiyi) + (zyitxiy_2)) \blacksquare x_2 - (ztxiyix_2) \blacksquare y_2$$
$$+ \{(x,y,txiyi)) + \{(x,ty,xiyi) - \{y_2txiyix_2)) \blacksquare z \tag{35}$$
$$+ (txiyix_2)Uy_{2tZ} - \{ty,xiyi)U_{X2tZ} = 0(IT\{n+3)),$$
$$(x,yitxiy_2) \blacksquare z + (zy,txiyi) \blacksquare x_2 + \{zx,tyixi) \bullet y_2 = 0(IT\{n+3)). \tag{36}$$

In this instance relations (31)-(34) are appropriate ^-representatives for $<_{71}$, $g_2$, $gs$, and $gg$, (35) and (36) are consequences of (31)-(34) modulo IT(n + 3).

Proof. (31) The definition of Shirshov-Cohn polynomials (see [2]) implies that

$$g1 =_2 (zx,y,x_2y_2) \blacksquare z - (z^2x,y,x_2y_2) - (xiyix,y_2z) \blacksquare z - \{zxiyix,y_2) \blacksquare z +$$
$$z^2x,y,x_2y_2) + z^2y_2x_2y,x_1)$$
$$= \{zxiyix,y_2) \blacksquare z - \{zy_2x_2yiXi) \blacksquare z - \verb|^|(Z^2XiyiX,y_2) + \verb|^|(Z^2y_2X,yiXi)$$
$$= (ZXiyiX_2y_2) \blacksquare z - \{zy_2x_2yiXi) \blacksquare z - \{z^2XiyiX_2) \blacksquare y_2 + (Z^2Xiyi)U_{x_2}y_2 - ]$$
$$\verb|^|\{x_2z^2xi(yiy_2)) + \{z^2y_2x_2yi) \blacksquare Xi - \{z^2y_2x_2)Uyi,xi + VIZ^2y2(x_2Xi)).$$

Therefore,

$$l\{gi) = \{zxiyix,y_2) \blacksquare t + (txiyix,y_2) \blacksquare z - \{zy,x,yix\{) \blacksquare t - (ty,x,yixi) \blacksquare z$$
$$\qquad (o$$
$$- \{\{ztxiyix_2) + (tzxiyix_2)) \blacksquare y_2 + ((zty,x,yi) + (tzy,x,yi)) \blacksquare X\verb|\|$$
$$= 0(1T \ (n+3)).$$

By using identities (30), we can readily ЈI-represent polynomials $\{zx\backslash y\backslash x,y_2)$-$t$ and $\{zy,x,y\backslash Xi) \blacksquare t$ as follows:

$$(zxiyix,y_2) \blacksquare t =_2 ((z \ o \ xi)yix,y_2) - \{xizyix,y_2) \blacksquare t$$
$$= 4(yix_2)U_{z \cdot Xli}y_2R_t - 2(zyix_2)U_{Xuy_2}Rt$$
$$' \ ^xi)^t,(y,x_2)Ry_2 \ \text{T} \ ^y^t,\{yi_{X_2}\}R_z - _xi \tag{37}$$
$$\verb|^|tU_{z \cdot xli}y_2 \ R\{y_{ix_2})_2X\backslash U_{(zy_{ix_2}^f;Ry_2}}$$
$$_2y_2U\{_{zyiX_2})_t tR_{xl} +$$
$$_2tU_{Xl}^y_jji(_{z}y_{lX2}y$$

Now

$$2tU_{xl},y, R(_{z}y \ 1Ж2) \qquad ``^{yiU_{zx2}} - \{xity2)$$
$$\widetilde{2}\{x\backslash ty,\}U_{ZtX_2}^R_{yi} \ _{(5)}^{yiRx_2} \ \text{T} \ ^{x\verb|^|\verb|^|}\{_xity_2),yi \ R_z$$
$$= (\{xity_2\}zyi \blacksquare x_2 + (\{xity_2\}x_2yi) \blacksquare z - (z\{xity_2\}x_2) \blacksquare yi$$

$$= (Xity_2x_2yi) \blacksquare z - (zxity_2x_2) \blacksquare m(W(n+3)), 4y_2 17u_{y1X2}) Rz\text{-}Xl$$
$$(0$$
$$— 2(f\text{Y}2\{f/lA2\}) \; (RzRxi \; 4 \cdot RxiRz \; Uxi \;_iz)$$
$$= ((zty_2X_2yi) + (ty_2x_2yiz)) \blacksquare X\backslash - \{(xity_2x_2yi)$$
$$(v$$
$$+ (ty_2X_2yiXi)) \blacksquare z - 2\{ty_2x_2yi) U_{Xl>z}(W\{n+3)).$$

Proceeding similarly with other terms in (37), we arrive at

$$(zxiyix_2y_2) \blacksquare t = (zty_2x_2yi) \blacksquare xi - 2(zxity_2x_2) \blacksquare yi + \{tzxiyix_2) \blacksquare y_2$$
$$+ (2(xity_2x_2yi) + (ty_2x_2yixi)) \blacksquare z \tag{38}$$
$$- 2\{ty_2x_2yi) U_{XuZ}(W(n+3)).$$

By substituting $x\backslash \;^\wedge y_2$ and $x_2 \;^\wedge y\backslash$ in (38), we obtain analogous A-representations for $\{zy_2x_2y\backslash Xi) \; \text{-} t$. If the resulting expressions are substituted in $l(gi)$ we are faced with relations (31). (32)-(34) Are verified similarly.

(35) We have

$$(\cdot txiyix_2) U_{y2Z} = 2(txiyi)R \; x_{,Uy2,z}$$
$$(0$$
$$=_{(3)} 2\{txiyi)(- \; \text{-}Ry_2U_x2,Z \; RzU_x2,Y2 \qquad Uy2_jzRx2 \; Uy_{2j}zRx2$$
$$U_{x2},zRy2 \; Ux2,y_2\sim^\wedge z)$$
$$=_{(0} \text{-}\{y_2txiyi) U_{X2)Z} \sim \{txiyiz) U_{X2j}y_2 + \{y_,txiyiz) \bullet x_2$$
$$+ \{\{x_2tx_1y_1z) + (ztxiyix_2)) \bullet y_2 + \{y_2txiyix_2) \blacksquare$$
$$z(W(n+3)).$$

Now

$$\{y_,tXiyi) U_{X2j}z — ^\wedge\{Y2^\wedge\} \; Rxi \; `yiU_{X2i}z \; \{ty_2Xiyi) U_{X2j}z \; (0$$
$$`` ^\wedge\{y2t)\{Uxi\text{-}yi_yX2\text{-}^\wedge Z \; Uxi\text{-}yi_yzRx2\sim) \; \{ty_2Xiyi) U_{X2},Z$$
$$К Ч$$
$$= \text{-}\{ty_2Xiyi) U_{X2j}z + (\{zty_2xiyi) + (zy_,txiyi)) \bullet x_2$$
$$\backslash 4$$
$$+ \{(x_,y_,txiyi) + \{x_,ty_,zxiyi)) \blacksquare z(W\{n+3)).$$

Analogously,

$$\{txiyiz)llx_2,y_2 \sim \;^\wedge\{^\wedge Xl) \; Ry_{,\cdot_2}Ux_2>y_2$$
$$К Ч$$
$$=_{4\cdot}(tXl)(Uy_{,\text{-}zj}x_,Ry_2 + Uyi'Z,y,Rx_2)$$
$$\textbf{(3)}$$
$$= (x_2txiyiz) \blacksquare y_2 + ((zyitxiy_2) + (zyiXity_2)) \blacksquare x_2(W(n+3)).$$

Therefore.

$$(ty_2xm)U_{X2>z} = ((zty_1xiyi) + (zy_1txiyi) + (zyitxiy_2)) \blacksquare x_2$$
$$- (ZtXiyiX_2) \blacksquare y_2 + \{(x_1y_1txiyi) + (x_1ty_1xiyi)$$
$$- (y_2tx_1y_1x_2) \blacksquare z + (txiyix_2)Uy_{2|Z}(W\{n+3)).$$

(36) From (35), we derive

$$(\blacksquare tx_1y_1x_2)Uy_{2|Z} = \{tyix_1y_2)U_{XuZ} + \{(ztxiyix_2) + (zxityix_2)$$
$$+ (zx_2tyixi)) \blacksquare y_2 - (ZtyiX_1y_2) \blacksquare X\backslash + \{\{y_1XityiX_2)$$
$$+ (y_2txiyix_2) - (xityix_1y_2)) \blacksquare z(W(n+3)).$$

Using (32) yields

$$(\blacksquare tx_1y_1x_2)Uy_{2|Z} = \{ty_1x_2yi)U_{XuZ} - \{zty_1x_2yi) \blacksquare X\backslash - \{zy_2txiyi)$$
$$\blacksquare x_2 + \{zxity_1x_2) \blacksquare y_2 + (ztxiyix_2) \blacksquare y_2 - ((xity_1x_2yi)$$
$$- (y_2tx_1y_1x_2)) \blacksquare z(W\{n+3)).$$

Therefore,

$$\{\{tyix_2y_2) - (ty_1x_2yi)U_{XuZ} + ((ztxiyix_2) \quad + (zxityix_2) + \{zx_2tyitxi)) \blacksquare y_2$$
$$- (ztyix_1y_2) \blacksquare Xi + \qquad\qquad \{\{y_1xityix_2) + (y_1txiyix_2)$$
$$- (xityix_1y_2)) \qquad \blacksquare \qquad\qquad\qquad z$$
$$+ (zty_1x_2yi) \blacksquare Xi + \{zy_1txiyi) \blacksquare x_2 - (ZXity_1X_2) \blacksquare yi$$
$$- \{ztxiyix_2) \blacksquare y_2 + \qquad\qquad ((xity_1x_2yi) \qquad - (y_2txiyix_2)) \blacksquare$$
$$z$$
$$= \{\{tyix_2y_2) - (ty_1x_2yi))U_{XuZ} + \{\{zty_1x_2yi) - (ztyix_1y_2)) \blacksquare Xi$$
$$+ (zy_2txiyi) \blacksquare x_2 - (ZXity_1X_2) \blacksquare yi + ((zxityix_2) + (zx_2tyixi))$$
$$\blacksquare y_2 + \{(xity_1x_2yi) - (xityix_1y_2) + (y_2xityix_2)) \blacksquare z = 0(W\{n+3)).$$

Relation (32) gives rise to

$$\{\{tyix_2y_2) - (ty_2x_2yi)U_{XuZ} + ((zty_2x_2yi) - (ztyix_2y_2)) \blacksquare Xi +$$
$$(zxityix_2) \blacksquare y_2 - (ZXity_2X_2) \blacksquare yi + ((xity_2x_2yi)$$
$$- (xityix_1y_2)) \blacksquare z = 0(W(n+3)).$$

Hence

$$(.zy_1txiyi) \blacksquare x_2 + (zx_2tyixi) \blacksquare y_2 + (y2Xityix_2) \blacksquare z = 0(W\{n+3)).$$

The proposition is proved.

ter for a derivational substitution $u_{-z}$. Using the $B_{-}$ representations above, we prove

LEMMA 6. For any $t = (\text{¾})$ with $\kappa(u_\ell) = n \wedge 1$, the following relations hold:

$$(zixityix_2) \; \blacksquare \; y_2 = 0(\mathrm{IF}(n+3)), \qquad (39)$$

$$(xityix_2y_2) \; \blacksquare \; zi = 0(\mathrm{IF}(n+3)), \qquad (40)$$

$$(xiyitx_2y_2) \; \blacksquare \; zi = 0(\mathrm{IF}(n+3)), \qquad (41)$$

$$(tyiXiy_2)U_{ZUX2} = (yitxiy_2)U_{ZUX2} = 0 \; (W(n+3)). \qquad (42)$$

Proof. If $h(ui) = 1$, then (39)-(42) have height 6 and follow from $G\%$ by the argument in item (6). Assume $I\mathcal{E}(u_\ell) \wedge 2$. Consider all possible types of $u_\ell$. In what follows, a reference to the definition of Shirshov-Cohn polynomials and the inductive assumption is marked with ($d$). If $\text{Щ} = xsViX_4$ or $y_3V_2y^\wedge$, $h(vi) \wedge 1$, then relations (39)-(42) follow from ($d$). There are three cases to consider:

(a) $u_\ell = z_2ViZ_3$, with $h\{vi\} \wedge 1$;

(b) $m = z_2ViX_3$ or $z_2Viy_3$, with $h(vi) \wedge \mathrm{o}$;

(c) $m = x_3V_2V_3$, with $h(vi) \wedge \mathrm{o}$.

Below, for brevity, we write $/ = g$ in place of $/ = g(W(n+3))$. First we argue for relations (39)-(41).

(a) Let $u_\ell = z_2ViZ_3$. From (32), we derive

$$\{zixity_2x\{\} \; \blacksquare \; yi - \{zix_2tyixi\} \; \blacksquare \; y_2 + \{\{x_2tyixiy_2\} - \{x_2ty_2xiyi\}\}) \; \blacksquare \; Zi$$
$$\underset{(32),\,(d)}{=} {}_2((x_2ty_2xiy_2) - (x_2tyixiyi)) \; \blacksquare \; Zi = \mathrm{o}$$

or

$$(x_2tyiXiy_2) \; \blacksquare \; Zi = (x_2ty_2xiyi) \; \blacksquare \; Zi.$$

Substituting $y\backslash = y_2 = y$ in (31) yields

$$(zytxiy) \; \blacksquare \; x_2 - (zxityx_2) \; \blacksquare \; y + ((xityxiy) - (ytxiyx_2))$$
$$\underset{(d)}{=} {}_2((xityx_2y) - (ytxiyx_2)) \; -z = \mathrm{o},$$

or

$$(xityx_2y) \; \blacksquare \; Zi = (ytxiyx_2) \; \blacksquare \; Zi.$$

Consequently, the element $(XityiX_2y_2) \; \blacksquare \; Zi$ is symmetric in $x\backslash$, $x_2$ and in $yi$, $y_2$ modulo $\mathrm{IF(n+3)}$, and

$$(XityiX_2y_2) \; \blacksquare \; Z = (yitXiy_2X_2) \; \blacksquare \; Zi.$$

Now (36) entails

$$(x_2yitxiy_2) \; \blacksquare \; z \underset{(d)}{=} -((y_2tXiyiX_2) \cdot (x_2tyiXiy_2)) \; \blacksquare \; z.$$

Hence the element $(o\char`^2/1\char`^22/2) \bullet z\backslash$ is symmetric in aq, $X_2$ and in $y_i$, $y_2$ modulo $W(n+3)$, and

$$(xiyitx_2y2) \;\blacksquare\; z_i = -2(x_ity_ix_2y2) \;\blacksquare\; z_i = \boldsymbol{-2(y_1tx_1y_2x_2)} \;\blacksquare\; z_i$$

$$= \boldsymbol{-2(y_1x_1y_2x_2)} \;\blacksquare\; z_i = \boldsymbol{-2\{xiyix_2ty2)} \;\blacksquare\; z_i \qquad (43)$$

$$= \boldsymbol{(yiXIty_2X2)} \;\blacksquare\; z_i.$$

Put $Xi = X_2 = x$ and $2/1 = 2/2 = 2/-$ Then

$$(xtyxy) \;\blacksquare\; z \underset{(d),(l)}{=} 4(ty)U_{xIx}.yR_z - (yxtyx)R_z$$

$$\underset{(Z),(d),(38)}{=} 4 \cdot (ty)(-U_{Z/X}.yR_x + R_x.yU_{x,z}) + 2(xtyxy) \;\blacksquare\; z$$

$$\underset{(l),(d)}{=} -(yxtyz) \;\blacksquare\; x + (yxtyy)U_{XtZ} + 2(xtyxy) \;\blacksquare\; z$$

$$\underset{(d)}{=} (yxxty) \;\blacksquare\; z + 2\ (xtyxy) \;\blacksquare\; z$$

$$(xtyxy) \;\blacksquare\; z + 2(xtyxy) \;\blacksquare\; z.$$

Consequently,

$$(xtyxy) \;\blacksquare\; z = o$$

$$(x,tyiX_2y2) \;\blacksquare\; Zi = (XiyitX_2y2) \;\blacksquare\; Z = \boldsymbol{0,}$$

$$(ZiXityiX_2) \;\blacksquare\; \boldsymbol{2/2} \underset{(d)}{=} -(xityiX_2y2) \;\blacksquare\; Zi = \boldsymbol{0.}$$

Case (a) is completed.

(b)   Let $u_t = Z_2ViXs$, with $h(vi)\char`^0$. We have

$$(ZityiXiy_2) - x_2 = (ZiX_3V^*Z_2yiXiy) \;\underset{кч}{\blacksquare}\; x_2 = 0,$$

$$(Ziy_2tXiyi) \;\blacksquare\; x_2 = (Ziy2X2V^*Z_2Xiyi) \;\underset{кч}{\blacksquare}\; x_2 = 0,$$

$$(ZiXity_,X_2) \;\blacksquare\; 2/1 \underset{(0)}{=} -(Xity_,X_2y2i) \;\blacksquare\; Zi, \qquad (44)$$

$$(zitXiyiX_2) \;\blacksquare\; 2/2 = (zix_,v^*z_,x_iy_ix_2) \;\underset{кч}{\blacksquare}\; 2/2 = {-}{<}\char`^1Ж_,г7^*\char`^_2Ж12/1Ж_3) \bullet 2/2-$$

Now if we apply (36) we obtain

$$(ZitXiyiX_2) \;\blacksquare\; 2/2 = (ZiX_3V^*Z_2XiyiX_2) \;\blacksquare\; 2/2 =$$

$$(ZiX_3\{v^*Z_2Xi\}yiX_2) \;\blacksquare\; 2/2$$

$$\underset{(36)}{=} -(x_3yitiX_2y2) \;\blacksquare\; Zi - (Ziy_2tiX2yi) \;\blacksquare\; x_3,$$

where $H = \{v^*Z_2Xi\}$. Therefore,

$$(ZitXiyiX_2) \;\blacksquare\; 2/2 = -(x_,yiXiZ_2ViX_2y2)\ 'И + {<}\char`^12/2Жi\char`^гЖ_,2/1)\ 'X_3$$

$$\underset{(0,(44)}{=} (X_2yiXiZ_2ViX_,y_2) \;\blacksquare\; Zi = (X_2yiXi\{Z_2ViX_,\}y_2) \;\blacksquare\; Zi \qquad (45)$$

$$= (X2yiXity_2) \;\blacksquare\; Zi.$$

By substituting the resulting relations in (31), (32), and (36), we arrive at

$$2((xity_2x_2yi) - (y_2txiyix_2)) \blacksquare zi + \{txiyix_2\}U_{y2>z} = 0, \tag{31}$$

$$\{x_2tyiXiy_2\} \blacksquare Zi = \{x_2ty_2xiyi\} \bullet Zi, \tag{32}$$

$$\{x_2yitxiy_2\} \bullet Zi = \{zx_2tyixi\} \bullet y_2 = (X_2tyiXiy_2) \bullet Zi - \tag{36} \tag{d}$$

Now we prove that

$$(xiz_2ViX3yix_2y_2) \blacksquare Zi = 0.$$

We have

$$(\blacksquare XlZ_2ViX_2yiX_2y2yi) \blacksquare Zi = 2\{z_2ViX;i\}x_1,y_1x_2yiRz$$

$$= 2\ (z_2ViX_3)(-\ Uzi,yixiyiRxi + U_{x1}\text{-}zi,yiX2yi + Ryix_2yiUxi,zi) \tag{5}$$

$$= \sim\{yiXiyiz_2ViX_2zi\} \blacksquare X\backslash + 4:\{z_2ViX_2\}U_{x1} \bullet Zl_yX_2\text{-}Vl\text{-}yi \\ K4$$

$$= \wedge\{z_2ViX_2\}(U_{Xl\cdot Zljx2}.y_1Ry_1 + U_{Xl\cdot Zlt}y_1R_{X2\cdot yi}) \tag{5),6)}$$

$$= \{ZiXiZ_2ViX_2yiX_2\} \blacksquare yi + \{x_2yiZ_2ViX_2ZiXi\} \blacksquare yi \tag{(i)}$$
$$+ \{\{ziXiz_2ViX_2yi) + (yiz_2v_1x_2zi,x_1))(R_{X2}Ry_1 + R_{y1}R_{x2} - U_{X2tV1})$$

$$= \frac{1}{(d),(l)} -\{XiZ_2ViX_2yiX_2yi) \blacksquare Zi - \{yiX_2yiZ_2ViX_2Zi) \blacksquare X\backslash \\ 2 \qquad\qquad\qquad l$$

$$+ \sim\{ZlXiZ_2ViX_3yiX_2)$$

$$\blacksquare yi =$$

Consequently,

$$(x_2tyiXiy_2) \blacksquare Zi = (x_2yitxiy_2) \blacksquare Zi = 0,$$

$$_2(y_2txiyix_2) \blacksquare Zi = (txiyix_2)Uy_{2>Zl}.$$

Ultimately,

$$(y_2tXiyiX_2) \blacksquare Zi = \{y_2X_3V^*Z_2XiyiX_2) \blacksquare Zi,$$

$$\{tXiyiX_2)U_{y2tZl} = (x_3v^*z_2XiyiX_2)Uy_{2tZl}.$$

If $Vi = WiZ_3$, then $(tXiyiX_2)U_{y2:Zl} = 0$, whence $(y_2tXiyiX_2) \blacksquare Zi = 0$. Hence $Vi = Wiy_3$. In this event
$$(o$$
$$(y_2tXiyiX_2) \blacksquare Zi = (ZitXiyiX_2) \blacksquare y_2 = (ZiX_3y_3W^*Z_2XiyiX_2) \ -1/2 = 0 \tag{45} \tag{(l)}$$

and $(tXiyiX_2)U_{y2tZl} = $ o.

(c) Let $u_i = x_3 y_r y_3$, with $h(v_i) \wedge 0$. We have

$$(x_2 y it Xiy_2) \blacksquare Zi = \sim(Ziy_2 tXiyi) \blacksquare X_2 - (ZiX_2 tyiXi) \, '2/2 = 0.$$
$$(36) \qquad\qquad (d)$$

Hence $(x_2 y it x iy_2) \blacksquare z \backslash = 0$. Now

$$(z it y ix_2, y_2) \blacksquare xi = \{ziy_3 V^* x_2 yix_2 y_2) \blacksquare Xi.$$

If $Vi = WiX_4$, then $\{z \backslash t y ix_2, y_2) \blacksquare x \backslash = 0$. If $Vi = WiZ_2$, then
$$(d)$$

$$(z it y ix_2, y_2) \blacksquare Xi = (z iy_3 \{v^* X_3 yi\} x_2, y_2) \blacksquare$$

$Xi = 0$ by the argument in case (b). Thus $\{z \backslash t y \backslash x_2, y_2) \blacksquare x \backslash = 0$.

Furthermore,

$$(x it y_2 x_2 yi) \blacksquare Zi = (xxy_3 V^* x_3 y_2 x_2 yx) \blacksquare Zx = \{x iy_3 \{v^* X_3 y_2\} x_2 yi) \blacksquare Zi$$
$$= \sim \{z iy i \{v^* X_3 y_2\} x_2 y_3) \blacksquare Xi - \{z iXi \{v^* X_3 y_2\} y_3 X_2)$$
$$\blacksquare yi \, (36)$$
$$= \sim \{z iy iv^* X_3 y_2 x_2 y_3) \blacksquare xi - (z_2 x_2 y_2 x_2 v_b y_3 X_2) \blacksquare y_x$$
$$= \sim \{z iy iv^* X_3 y_2 x_2 y_3) \bullet x_b \, (d)$$

If $Vi = WiX_4$, then $\{z \backslash y \backslash V^* x_3 y_2 x_2 y_3) \blacksquare x \backslash = 0$. If $Vi = WiZ_2$, then $\{z iy i \{v^* X_3 y_2\} x_2 y_3) \blacksquare x \backslash = 0$
by
$$(d)$$

the argument in case (b).

    Relations (39)-(41) are thus proved.

    (42) In view of (39)-(41), we have

$$\{t x iy ix_2) Uy_{2tZ} \underset{(31)}{=} \{t y ix_2, y_2) U_{xuz}.$$

Therefore, $\{t x \backslash y \backslash x_2) U_{y_2 z_1} = 0$ if $u_i = z_2 ViZ_3$, $z_2 ViX_3$, $z_2 Viy_3$.

    Consider the case where $u_i = X_3 V_2 y_3$, with $h(vi) \wedge 0$. Now (32), combined with (39)-(41), implies $\{t x \backslash y ix_2) U_{v2:Zl} = \{t x_2 y \backslash Xi) U_{v2:Zl}$. Consequently, the element $(t x \backslash y ix_2) U_{v2Zl}$ is symmetric in $X \backslash$, $x_2$ modulo $IT(n+3)$. Consider

$$(t x iy iXi) U_{y2 > z l \backslash} \underset{4}{\vert} 2 \{t) Rxiyixi \, U_{y_2, z \backslash}$$

$$\underset{(5)}{=} )(\wedge 2/2 , XiyiXlRzi \, T \, U_{z_1} \, _{t} xiyiXl \, Ry_2 \qquad Uxiyixi, y_2 - Zl ,$$

$$= 4(1)(7 \, y_r z \backslash, x \backslash - y \backslash - x \backslash$$
$$(39)-(41)$$

$$= \sim 4(t)(C_{y_r z l \backslash, x \backslash - y \backslash R - x \backslash} + U_{y_r z_2 i, xiRxi - yi)}$$
$$(0.(5)$$

where $x\backslash y\backslash X\backslash = 2x\backslash \blacksquare y\backslash \blacksquare x\backslash — x\backslash \blacksquare y\backslash$. Hence $(txiyix_2)U_{y_2tZ} = 0$ for all $u_t$. Finally,

$$\langle Xltyix_i\rangle Uy_{2,zi} - 4(\divideontimes 17) R yi\text{-}X2 \,\rangle \,{}^{\wedge} U2, zi - (\%2yi\%lt) U y_{2,zi}$$

$$= 4\,(x,t)\,(U_{yi\text{-}X_2\text{-}X_2\,R}\dots)\dots)$$

$$\underset{(39)\text{-}(41)}{\equiv} -\langle x_2 y_1 x_1 t\rangle U_{y_2,z} \equiv 0.$$

The lemma is completed.

**LEMMA 7.** For any $t = (\frac{3}{4})$ with $h(ui) \wedge 1$, the following relations hold:

$$(xitx_2yizi) \blacksquare y_2 = 0(W(n +$$
$$3)), (x_2tyixizi) \blacksquare y_2 = 0(W(n$$

**Proof.** As in the proof of Lemma 6, we have $h(ui) \wedge 2$. Consider the following cases

(a)   $Щ_1 = z_2 ViZ_3]$

(b)   $m = z_2 ViX_3$ or $m = z_2 Viy_3;$

(c)   $u_t = X_3 Viy_3.$

The cases with $u_t = x_3 ViX_4$ or $y_3 Viy_4$ are trivial in view of (34).

(a) Let $u_t = z_2 ViZ_3$. Then

$$\{XiZ_2 ViZ_3 X_2 yiZi) \blacksquare y_2 = \{xiz_2\{v*Z_3 X_2\}yiZi) \blacksquare y_2 - \{XiZ_2 X_2 Z_3 V*yiZi) \blacksquare y_2$$
$$= \sim \{XiZ_2 X_2 Z_3 V* yiZi) -y_2$$
$$\underset{\text{(оУ) (41J}}{}$$
$$= \sim \{xiz_2 x_2\{z_3 v*yi\}zi) \blacksquare y_2 \underset{(39)\text{-}(41)}{=} 0.$$

Consequently,

$$(xitx_2yiZi)\; 2/2 = 0.$$

Furthermore.

$$('tx_2yiZi)U_m\,,Xl = (tx_2yi)R\; ZlUy_2,Xl$$
$$\underset{3)}{=} 2\langle tx_2 y_1\rangle(-R_{y_2}U_{z_1,x_1} - R_{x_1}U_{y_2,z_1} + U_{y_2,x_1}R_{z_1}$$
$$+ U_{z_1,x_1}R_{y_2} + U_{y_2,z_1}R_{x_1})$$
$$\underset{(1),(39)\text{-}(42)}{,} = (xitx_2y_{1Zi})Ry_2 + (y_2tx_2yiZi)R_{xl}\underset{(48}{=} \{y_2tx_2yiZi)R\; XI\;*$$

On the other hand

$$(z_2v_iz_3x_2y_1z_1)Uy_{2\}Xl} = (z_2\{v_iz_3x_2\}y_1z_1)Uy_{2\}Xl} - \{z_2x_2z_3v*yiZi)U_{y2tXl}\underset{(0,(42)}{}\; 0.$$

Therefore, (46) and (47) are satisfied in case (a).

(b) Let $u_4 = z_2 v i x_3$. Then

$$(x i t x_2 y i z i) \; \mathbf{2}/2 = 0,$$

$$(z_1 y_1 x_2 t y_2) \; \blacksquare \; x i = (z i y i x_2 z_2 V i X_3 y_2) \, -x i = 0,$$

$$\{x_2 t y_2 x i z i) \; \blacksquare \; y i = 4 : (t y_2) U_{x_2} \, , x_1 -$$
$$z_1 R y_1 \; _{(0}$$
$$= \; _{U M^5)} \; {}^{\wedge} \{ t y_2) \, ( \; - U y_1 , x_1 - z_1 R x 2 + R x_1 - z_1 U_{X_2}, y i)$$
$$_{(42)} = -[V l t V 2 X l Z i ) R_{X2} + (t y_2 X i Z i) U_{X 2 m}. \; _{(/), (39)-}$$

On the other hand,

$$(. Z_2 V i X_3 y_2 X i Z i) U_{X 2 > y i} = -(z_2 y_2 x_3 v * x i Z i) U_{X 2 > y i} = 0.$$
$$_{(34)} \qquad \qquad _{(/)}$$

Now (34) implies that $(x_2 t y_2 x \backslash Z i) \; \blacksquare \; y \backslash = \{y i t y_2 x \backslash Z \backslash) \; \blacksquare \; x_2 = 0$.

(a) Let $u_4 = x s V i y s$. From (34), we derive

$$(y 2 t x_2 y_1 z_1) \; \blacksquare \; x i + (x i t y i x_2 z i) \; * 2/2 = o.$$

Since

$$\{y 2 t X_2 y i Z i) \; \blacksquare \; X i = [y 2 X_3 V i y_3 X_2 y \backslash Z i) \; \blacksquare \; X \backslash = 0,$$

we have $\{x \backslash t y \backslash x_2 z \backslash) \; \blacksquare \; y_2 = 0$. The lemma is completed.

## 7. BASIC RESULT

T-ideals $W(n), \; n \wedge 2$, are described in

**LEMMA 8.** For any $(\tfrac{3}{4}) \; G \; R$ with $h(u i) = n \wedge 6$, we have

$$2 (\tfrac{3}{4}) \tfrac{3}{4} - (o \tfrac{3}{4}) - (\text{«ja}) \; G \; IF(n - 1), \qquad 2(u_4) u_{a > b} - (a u_4 b) - (b u_4 a) \; G \; W(n). \qquad (49)$$

**Proof.** By virtue of Lemmas 4-7, it suffices to verify the following relations:

$$A = (z_1 y_1 z_2 u_1 x_1 y_3) \; \blacksquare \; x_2 = 0(W\{n + 4)),$$
$$B = \{x i y i U i X_2 y_2) \; \blacksquare \; Z i = 0(W\{n + 3)),$$
$$C = (x i y i U i y_2 x_2) U_{z i m} = 0 \; (W \; (n + 4)),$$

where $i \partial (u_4) \wedge o$.

If $h(u i) = 0$, then relations $A$-$C$ have height 6 and follow from $Gg$ (see Sec. 5). Let $I \partial (u_4) \wedge I$. We verify $A$. If $u_4 = V i y_3$, then

$$A = \{z i y i z_2 V i y ; i X i y_2) \; \blacksquare \; x_2 = (z i y i \{z_2 V i y_3\} x i y_2) \; \blacksquare \; x_2 = \qquad 0(W(n + 4)).$$
$$_{(o \tilde{y}) - (41)}$$

If $u_i = v i z_3$, then

$$A = \{z i y i Z_2 V i Z_3 X i y_2\} \blacksquare x_2 = (z_1\{y_1 Z_2 V_1\} z_3 X_1 y_2) \blacksquare x_2 - (Z i V^* Z_2 y i Z_3 X i y_2) \blacksquare x_2$$

$$= -(z i v^* z_2 y_4 z_3 x_4 y_2) - x_2\{W \{n+4\}).$$

**(46)**

**If** $v^* = x_3 W i, y_3 x_3 W i$, **then** $A = \mathbf{0}(W\{n+4\})$. **The only case left is** $v^* = y_3 z_4 W i$. **In this instance** **(0**

$$A = -(z_1 y_3 z_4 w_1 z_2 y_1 z_3 x_1 y_2) \blacksquare x_2$$

$$= -\{z m\{z_4 W i Z_2 y i z_3\} x i y_2) \blacksquare x_2 +$$

$$\{z i y_3 z_3 y i z_2 w^* z_4 x i y_2) \blacksquare x_2 = \sim\{z i y_3 \{z_4 W i Z_2 y i Z i\} x i y_2)$$

$$\blacksquare x_2 + (z_1 y_3 z_3\{y_1 z_2 w^* z_4 x_1\} y_2) \blacksquare x_2 = \qquad 0(W\ (n+4)).$$

**(39)-(41)**

We argue for $B$. **If** $u_i$ **equals** $V i y_3$ **or** $x_3 V i$, **then**

$$B = (x_4 y i\{u i\} x_2 y_2) \blacksquare Z i = \mathbf{0}(\mathbf{IT}(n+3)).$$

**If** $Vi = z_2$, **then** $B = 0(W(n+3))$, **since** $h(B) = 6$. **It remains to handle the case where** $Vi = z_2 V i Z_3$, **with** $h(vi) \wedge 1$. **If** $Vi = W i X_3$, **where** $Wi$ **may be empty, then**

$$B = (x i\{y i z_2 W i X_3\} z_3 x_2 y i) \blacksquare z_4 = \qquad \mathbf{0}(\mathbf{IT}(n+3)).$$

**(46),(47)**

**Therefore, either** $Vi = y_3$, **or** $Vi = y_3 W i y_4$, **where** $h(wi) \wedge 1$. **In either case,**

$$B = \{x i y i z_2\{y_3 W i y_4 z_3 x_2\} y_2\} \blacksquare z_4 = \mathbf{0}(\mathbf{IT}(n+3)).$$

**(46),(47)**

We verify $C$. **In this event**

$$C = \{x i y i\{u i y_2\} x_2) U_{Z u y 3} = 0(W\ (n+4)).$$

**(42)**

**The lemma is completed.**

**THEOREM 4. All Jordan s-identities in three variables are consequences of the Glennie identity** $Gg$ **of degree 8.**

**Proof. Relations (48) imply that for** $n \wedge 6$,

$$J\ R(n)\ \mathbf{C}\ R(n-1) + B(n-2) + K,$$

$$[U(n)\ \mathbf{C}\ R(_n) + U(n-1) + K.$$

**Therefore,**

$$W(n) = R(n) + U(n-1) + K$$

$$\mathbf{C}\ R(n-1) + U(n-2) + R(n-1) + U(n-2) + K$$

**C** $W(n - 1)$.

Hence $W(n)$ C IF(5) + $K$ for any $n \wedge 6$. By the argument in Sec. 5, $W\{4) = 0$ and $W\{5)$ C $T(Gs)$- Therefore, $W(n)$ C K for all $n$ G N. In view of Theorem 2, $K = T(G\%)$. Hence $W(n)$ C $T(Gs)$. Now

$$\mathbf{T(Gs) \ C5_3 = \underset{n=2}{\overset{\infty}{J}|IT(n)} \ CT(G_8).}$$

Consequently, N$_3$ = $T(Gs)$. The theorem is proved.

Acknowledgments. I am grateful to Profs. I. P. Shestakov and K. McCrimmon for useful pieces of advise and helpful discussions.